\documentclass[12pt,draft]{article}

\usepackage{amsthm,amssymb,latexsym,amsmath}

\newtheorem{theorem}{Theorem}[section]

\newtheorem{proposition}[theorem]{Proposition}
\newtheorem*{remark}{{\it Remark}}

\newcommand{\nc}{\newcommand}
\nc{\N}{{\mathbb N}}
\nc{\Z}{{\mathbb Z}}
\nc{\R}{{\mathbb R}}
\nc{\C}{{\mathbb C}}
\nc{\HH}{{\mathbb H}}
\nc{\dd}{{\rm d}}
\nc{\slc}{{\mathfrak s}{\mathfrak l}(2,\C )}
\nc{\su}{{\mathfrak s}{\mathfrak u}(2)}
\nc{\slr}{{\mathfrak s}{\mathfrak l}(2, \R )}
\nc{\g}{{\mathfrak g}}
\nc{\cg}{{\mathfrak g}^{\C}}

\begin{document}

\title{Gravitational interpretation of the Hitchin equations}
\author{G\'abor Etesi\\
\small{{\it Department of Geometry, Mathematical Institute, Faculty of 
Science,}}\\ 
\small{{\it Budapest University of Technology and Economics,}}\\ 
\small{{\it Egry J. u. 1, H \'ep., H-1111 Budapest, Hungary
\footnote{{\tt etesi@math.bme.hu}}}}\\
\small{ and}\\
\small{{\it Instituto de Matem\'atica, Estat\'\i stica e 
Computa\c{c}\~ao Cient\'\i fica,}}\\
\small{{\it Universidade Estadual de Campinas,}}\\ 
\small{{\it C.P. 6065, 13083-859, Campinas, SP, Brazil\footnote{{\tt 
etesi@ime.unicamp.br}}}}}

\maketitle

\pagestyle{myheadings}
\markright{G. Etesi: Gravitational interpretation of the Hitchin 
equations}

\thispagestyle{empty}

\begin{abstract}
By referring to theorems of Donaldson and Hitchin, we exhibit a rigorous 
AdS/CFT-type correspondence between classical $2+1$ 
dimensional vacuum general relativity theory on $\Sigma\times\R$ and 
SO$(3)$ Hitchin theory (regarded as a classical conformal field theory) on 
the spacelike past boundary $\Sigma$, a compact, oriented 
Riemann surface of genus greater than one. Within this framework we can 
interpret the $2+1$ dimensional vacuum Einstein equation as a decoupled 
``dual'' version of the 2 dimensional SO$(3)$ Hitchin equations. 

More precisely, we prove that if over $\Sigma$ with a fixed conformal 
class a real solution of the SO$(3)$ Hitchin equations with induced flat 
SO$(2,1)$ connection is given, then there exists a certain 
cohomology class of non-isometric, singular, flat Lorentzian metrics 
on $\Sigma\times\R$ whose Levi--Civita connections are precisely the lifts 
of this induced flat connection and the conformal class induced by 
this cohomology class on $\Sigma$ agrees with the fixed one.

Conversely, given a singular, flat Lorentzian metric on $\Sigma\times\R$ 
the restriction of its Levi--Civita connection gives rise to a real 
solution of the SO$(3)$ Hitchin equations on $\Sigma$ with respect to the 
conformal class induced by the corresponding cohomology class of the 
Lorentzian metric.
\end{abstract}

\centerline{\small{AMS Classification: Primary: 53C50; Secondary: 58J10, 
83E99}}

\section{Introduction}

The aim of this paper is to offer a new physical interpretation of the 
2 dimensional SO$(3)$ Hitchin equations as the coupled ``dual'' version of 
the $2+1$ dimensional vacuum Einstein equation. This interpretation 
emerges within a rigorous AdS/CFT-type correspondence between Hitchin's 
theory of Higgs bundles over Riemann surfaces, regarded as a $2$ 
dimensional conformal field theory on the boundary and classical 
$2+1$ dimensional vacuum general relativity on a bulk space-time. Our work 
may be viewed as an arch spanned by Witten's ideas \cite{wit1} between two 
massive bearers: a theorem of Hitchin \cite{hit} and another one by 
Donaldson \cite{don} respectively, as follows.

The relationship between three dimensional general relativity and
gauge theory is not new. The main motivation for studying their 
common features comes from the effort to formulate a satisfactory quantum 
theory underlying three dimensional general relativity (we cannot make 
the attempt to survey the vast literature of the quantization issue; 
for a survey cf. \cite{lol} or a more recent excellent one is \cite{car}). 
Just to mention one trial, Witten argued that Lorentzian vacuum general 
relativity theory should be equivalent to an ISO$(2,1)$ Chern--Simons theory 
at the full quantum level; thereby general relativity in $2+1$ dimensions 
turned out to be not only exactly soluble at the classical level but also 
renormalizable as a quantum field theory \cite{wit1}\cite{wit3}. The key 
technical tool in Witten's approach is formulating general relativity in 
terms of a connection and a ``dreibein'', instead of a metric.
This approach is remarkable because establishing any conventional 
relationship between general relativity and Yang--Mills theory apparently 
fails in any other dimensions, despite the efforts made over the past 
thirty years. 

More recently there has been interest among physicists in 
understanding the celebrated Maldacena conjecture or AdS/CFT 
correspondence \cite{mal} which sheds new light onto the gauge 
theory-gravity duality. Broadly speaking, this conjecture states the 
existence of a duality equivalence between some quantum gravitational 
theories on an anti-de Sitter space $M$ and quantum conformal field 
theories on the boundary at conformal infinity $\partial M$. At the 
semi-classical level and using a (Wick rotated) pure gravity theory in the 
bulk the correspondence was formulated by Witten \cite{wit2} and states that
\begin{equation}
Z_{CFT}([\gamma ])=\sum {\rm e}^{-I(g)}
\label{particiofv}
\end{equation}
where $Z_{CFT}$ is the partition function of some conformal field theory 
attached to a conformal structure $[\gamma ]$ on $\partial M$ and $I$ is 
the (renormalized) Einstein--Hilbert action of an Einstein 
metric on $M$ with conformal infinity $[\gamma ]$. The formal sum is 
taken over all manifolds and metrics $(M,g)$ with given boundary data 
$(\partial M, [\gamma ])$. As we mentioned, three 
dimensionality is distinguished in the conventional understanding of 
the gauge theory-gravity relationship \cite{wit1}, in the holographic 
approach however, Chern--Simons theory may play a role in various dimensions 
\cite{bel-moo}.

The AdS/CFT correspondence at least in its strict classical form has 
attracted some attention from the mathematician's side as well and led to nice 
geometrical results (cf. \cite{and} for a survey and 
references therein). This paper can also be regarded as an attempt to 
extend further its mathematical understanding by a natural 
generalization of the very core of the correspondence as well as link 
three dimensional gravity with Hitchin's theory of Higgs bundles over 
Riemannian surfaces. Since this is an integrable system, the relationship 
provides a further explanation, different from Witten's, why three 
dimensional gravity is exactly soluble.

This link is probably not surprising because in our opinion three dimensional 
{\it classical} gravity in {\it vacuum} is a two rather than three 
dimensional theory in its nature as can be seen by a simple topological 
argument. In three dimensions a Ricci flat space is flat. Although every 
compact, orientable three-manifold has zero Euler characteristic i.e., 
admits Lorentzian structures, only six of them are flat: These are the six 
orientable compact flat three-manifolds and are not interesting examples 
of three-manifolds because all of them are just finitely covered by the 
three-torus. Consequently we have to seek non-compact spaces carrying 
solutions of the classical vacuum Einstein equation like the annulus 
$\Sigma\times\R$; however this is rather a two dimensional object from a 
topological viewpoint. Nevertheless three dimensionality enters at the 
full {\it quantum} level as it was pointed out by Witten \cite{wit3}. 

Our paper is organized as follows. In Section 2 we prove that a 
real solution of the SO$(3)$ Hitchin equations over a compact 
oriented Riemann surface $\Sigma$ of genus $g>1$ induces a certain 
cohomology class of singular solutions of the Lorentzian 
vacuum Einstein equation over the annulus $\Sigma\times\R$. 
``Singular'' in this context means that the metrics may degenerate. Such 
solutions appear naturally if the metric is expressed via a connection and 
an independent dreibein, as was suggested by Witten \cite{wit3}. Our 
construction is based on a theorem of Hitchin \cite{hit} 
(cf. Theorem \ref{hitchintetel} here) on the relationship 
between real SO$(3)$ Hitchin pairs and flat SO$(2,1)$ connections.

Conversely, in Section 3, by referring to a theorem of Donaldson 
\cite{don} (cf. Theorem \ref{donaldsontetel} here) which states the 
equivalence between flat PSL$(2,\C )$ connections and SO$(3)$ Hitchin pairs, we 
present the reversed construction, namely, starting from a cohomology 
class of flat singular Lorentzian metrics on $\Sigma\times\R$ one can 
recover a unique real solution of the SO$(3)$ Hitchin equations on $\Sigma$, 
regarded as the past boundary of the annulus. 

We present our main results in Section 4. On the one hand we exhibit the 
correspondence in a precise form using the field equations (cf. Theorem 
\ref{ads-cft}). This presentation has the remarkable feature that it 
exhibits the vacuum Einstein equation as a sort of ``untwisted'', or
``decoupled'' variant of the Hitchin equations and vica versa 
if we interpret the flat SO$(2,1)$ connection on the gravitational side 
as being dual to the (non-flat) SO$(3)$ connection on the gauge 
theoretic side and simlarly, we regard the dreibein as being dual to the 
Higgs field.

Then we rephrase the correspondence in terms of the metric on the 
gravitational side. This way we can see that this correspondence looks like a 
generalized geometric AdS/CFT correspondence. Here ``generalized'' 
means that the conformal geometry on the boundary emerges in a rather 
abstract way compared with the usual AdS/CFT correspondence, namely, by 
exploiting the conformal properties of an equation which is 
essentially the massless Dirac equation as explained by Propositions 
\ref{megoldasok} and \ref{teichmuller}. 

Finally, we conclude with some speculations on the subject in Section 5.

\section{From Hitchin pairs to flat metrics}

The embedding SL$(2,\R )\subset{\rm SL}(2,\C )$ induces the 
factorized embedding ${\rm PSL}(2,\R 
)\subset{\rm PSL}(2,\C )$ and we will write ${\rm PSL}(2,\R )\cong {\rm 
SO}(2,1)$. Let $(\Sigma , [\gamma ])$ be a 
compact, oriented Riemann surface of genus $g>1$ with the conformal 
equivalence class of a smooth Riemannian metric 
$\gamma$ that is, a complex structure on it. Moreover let $P$ be an 
SO$(3)$ principal bundle over 
$\Sigma$ with either $w_2(P)=0$ or $w_2(P)=1$ and denote by $P^\C$ the
corresponding complexified PSL$(2, \C )$ principal bundle. Regarding 
SO$(3)$ as a real form of PSL$(2,\C )$ there is an 
associated anti-involution $*$ of the complex Lie algebra $\slc$. 
If $\nabla_A$ is an SO$(3)$ connection with curvature $F_A$ on $P$ and 
$\Phi\in\Omega^{1,0}(\Sigma , {\rm ad}(P^\C))$ is a complex Higgs field then 
the Hitchin equations over $(\Sigma ,[\gamma ])$ read as follows 
\cite{hit}:  
\begin{equation}
\left\{\begin{array}{ll}                
                 F_A+[\Phi , \Phi^*] & =0 \\ [2mm]
                 \overline{\partial}_A\Phi & =0. \\ [2mm]
\end{array}\right.
\label{hitchin}
\end{equation}
Recall that these equations are the dimensional reduction of the four 
dimensional SO$(3)$ self-duality equations hence are 
conformally invariant and exactly soluble.

Consider a solution $(\nabla_A, \Phi )$ of (\ref{hitchin}) associated to 
a fixed SO$(3)$ principal bundle $P$. If ${\cal A}(P)$ is the affine space 
of SO$(3)$ connections 
over $P$ then a map $\alpha :{\cal 
A}(P)\times\Omega^{1,0}(\Sigma , {\rm ad}(P^\C ))\rightarrow 
{\cal A}(P^\C )$ is defined as 
\begin{equation}
\alpha (\nabla_A, \Phi )=\nabla_A+\Phi +\Phi^*.
\label{lapos}
\end{equation}
Clearly the map descends to the gauge equivalence 
classes. Locally, on an open subset the resulting PSL$(2,\C )$ connection 
$\nabla_B$ looks like $\nabla_B\vert_U=\dd +B_U$ with $B_U=A_U+\Phi 
+\Phi^*$. It is easy to see that $\nabla_B$ is flat. Indeed, one 
quickly calculates
\[ F_B=\dd B+B\wedge B= F_A+[\Phi 
,\Phi^*]+\overline{\partial}_A\Phi +\partial_A\Phi^*=0\]
via (\ref{hitchin}). 
One may raise the question of $\nabla_B$ is moreover 
real valued i.e., wether takes its value in SO$(2,1)$. This is answered by 
a theorem of Hitchin (cf. Proposition 10.2 and 
Theorem 10.8 in \cite{hit}) which is the starting point of our discussion:

\begin{theorem}{\rm (Hitchin, 1987)} Let $(\Sigma ,[\gamma ])$ 
be a compact, oriented Riemann surface of genus $g>1$ endowed with a 
conformal equivalence class $[\gamma ]$ of a smooth Riemannian 
metric $\gamma$. Let $P$ be a principal SO$(3)$ bundle over $\Sigma$ 
satisfying either $w_2(P)=0$ or $w_2(P)=1$. Denote by ${\cal M}(P)$ the 
moduli space of gauge equivalence classes of smooth solutions to the 
SO$(3)$ Hitchin equations over $P$ with respect to $[\gamma ]$. Consider a 
map $\sigma :{\cal M}(P)\rightarrow{\cal M}(P)$ given by
\[\sigma ([(\nabla_A, \Phi )]):=[(\nabla_A, -\Phi )].\]
The fixed point set of $\sigma$ has connected components ${\cal 
M}_0$, ${\cal M}_2$, ${\cal M}_4,\dots ,{\cal M}_{2g-2}$ for $w_2(P)=0$ 
and ${\cal M}_0$, ${\cal M}_1$, ${\cal M}_3,\dots ,{\cal M}_{2g-3}$ for 
$w_2(P)=1$. The subset ${\cal M}_0$ is the space 
of flat SO$(3)$ connections on $P$ while ${\cal M}_k$ with $k>0$ can be 
identified with the $6g-6$ dimensional moduli of smooth, flat, irreducible 
SO$(2,1)$ connections of the form (\ref{lapos}) on certain principal 
SO$(2,1)$ bundles $Q_k$ of Euler class $k$ over $\Sigma$. 
$\Diamond$
\label{hitchintetel}
\end{theorem}
\begin{remark}\rm
Putting a suitable complex structure $J$ onto ${\cal M}(P)$ 
the map $\sigma$ can be regarded as an anti-holomorphic involution i.e., a 
{\it real structure} on $({\cal M}(P),J)$. This explains why 
complex flat connections of the form (\ref{lapos}) corresponding to the 
fixed point set of $\sigma$ inherit a real nature in the sense above (cf. 
\cite{hit}, Section 10 for details). Notice that {\it all} irreducible, 
flat SO$(2,1)$ connections with non-zero Euler class over $\Sigma$ arise 
this way.
\end{remark}
 
\noindent We wish to use these real, flat, irreducible connections to 
construct 
certain flat Lorentzian metrics over $\Sigma\times\R$ with a fixed 
orientation induced by the orientation of $\Sigma$. Consider the 
standard 3 dimensional real, irreducible representation $\varrho :{\rm 
SO}(2,1)\times\R^3\rightarrow\R^3$ and take the associated real rank 3 
vector bundles $E_k:=Q_k\times_\varrho \R^3$.
We restrict attention to the bundle $E_{2g-2}$ for which we 
have an isomorphism $E_{2g-2}\cong T\Sigma\oplus\underline{\R}$. For 
simplicity we shall denote this bundle as $E$ and the associated 
irreducible, flat SO$(2,1)$ connections of Theorem \ref{hitchintetel} on $E$ 
as $\nabla_B$. 

Let $\pi :\Sigma\times\R\rightarrow\Sigma$ be the obvious 
projection and consider the pullback bundle $\pi^*E$. This bundle admits 
an irreducible, flat SO$(2,1)$ connection $\pi^*\nabla_B$ by pulling back 
$\nabla_B$ from $E$. In the remaining part of the paper we will study the 
complex
\begin{equation}
0\rightarrow\Omega^0(\Sigma\times\R
,\:\pi^*E)\stackrel{\pi^*\nabla_B}{\longrightarrow
}\Omega^1(\Sigma\times\R
,\:\pi^*E)\stackrel{\pi^*\nabla '_B}{\longrightarrow
}\Omega^2(\Sigma\times\R ,\:\pi^*E)\rightarrow 0
\label{felhuzas}
\end{equation}
where $\nabla '_B$ is the induced connection, and in particular its first 
cohomology
\[H^1(\pi^*\nabla_B)=\frac{{\rm Ker}\:(\pi^*\nabla '_B)}{{\rm 
Im}\:(\pi^*\nabla_B)}\]
which will turn out to be of central importance to us.

Let $\widehat{\pi^*E}$ be an affine vector bundle over 
$\Sigma\times\R$ whose underlying vector bundle is $\pi^*E$. Fix an 
element $\xi\in \Omega^1(\Sigma\times\R, \:\pi^*E)$. We can regard $\xi$ 
as a fiberwise translation in $\widehat{\pi^*E}$ by writing 
$\xi_X(\hat{v}):=\hat{v}+\xi (X)$ with $X$ a vector field on 
$\Sigma\times\R$ and $\hat{v}\in\Omega^0(\Sigma\times\R ,\: 
\widehat{\pi^*E})$. Therefore we have an embedding 
\[\Omega^1(\Sigma\times\R, 
\:\pi^*E)\subset\Omega^1(\Sigma\times\R ,\:{\rm End}(\widehat{\pi^*E})).\]
Consequently if $\pi^*\nabla_B$ is a flat SO$(2,1)$ connection on $\pi^*E$ 
then 
\begin{equation}
\hat{\nabla}_{B,\xi}:=\pi^*\nabla_B+\xi
\label{isolapos}
\end{equation}
is an ISO$(2,1)$ connection on $\widehat{\pi^*E}$ where ISO$(2,1)$ denotes 
the $2+1$ dimensional Poincar\'e group. Its curvature is
\[\hat{F}_{B,\xi}=\pi^*F_B+(\pi^*\nabla '_B)\xi +\xi\wedge\xi 
=(\pi^*\nabla '_B)\xi\]
since translations commute. We obtain that $\hat{\nabla}_{B,\xi}$ is flat 
if and only if
\begin{equation}
(\pi^*\nabla '_B)\xi =0.
\label{cartan}
\end{equation}
Observe that {\it all} irreducible, flat ISO$(2,1)$ connections arise this 
way. Clearly the gauge equivalence class of $\hat{\nabla}_{B,\xi}$ is 
unchanged if $\pi^*\nabla_B$ is replaced by an SO$(2,1)$ gauge equivalent connection 
and $\xi$ by $\xi +(\pi^*\nabla_B)v$ with an arbitrary section 
$v\in\Omega^0(\Sigma\times\R ,\:\pi^*E)$. Therefore on the one hand we 
identify the space $H^1(\pi^*\nabla_B)$ with the underlying vector space 
of gauge equivalence classes of smooth flat ISO$(2,1)$ connections on 
$\widehat{\pi^*E}$ with fixed SO$(2,1)$ part $\pi^*\nabla_B$.

On the other hand out of a flat ISO$(2,1)$ connection one can construct a 
singular Lorentzian metric on $\Sigma\times\R$ as follows \cite{wit1}. Fix 
once and for all a smooth SO$(2,1)$ metric $h$ on $\pi^*E$ 
(notice that this bundle is an SO$(2,1)$ bundle) and pick up a flat 
connection of the form (\ref{isolapos}) on $\widehat{\pi^*E}$. Via the 
isomorphism
\[\Omega^1(\Sigma\times\R, \:\pi^*E)\cong\Gamma ( {\rm 
Hom}(T(\Sigma\times\R ), \:\pi^*E))\]
we can interpret $\xi$ as a ``dreibein'' $\xi :T(\Sigma\times\R 
)\rightarrow\pi^*E$ (since $\pi^*E$ 
and $T(\Sigma\times\R )$ are isomorphic bundles). Assume for a moment that 
$\xi_x^{-1}:(\pi^*E)_x \rightarrow 
T_x(\Sigma\times\R )$ exists for all $x\in\Sigma\times\R$ that is, 
$\xi$ is invertible as a bundle map. Using $\xi$ we 
can construct a smooth Lorentzian 
metric $g_\xi :=h\circ (\xi\times\xi )$ on $T(\Sigma\times\R )$. 
Locally $(g_\xi )_{ij}=\xi_i^p\xi_j^qh_{pq}$. We can suppose that the 
metric constructed this way is inextensible. The connection $\xi^{-1}\circ 
(\pi^*\nabla_B)\circ\xi$ is compatible with $g_\xi$ hence it represents 
the Levi--Civita connection of $g_\xi$ if it is torsion free. However this 
is provided by (\ref{cartan}) since this equation is just the Cartan equation 
for the smooth metric $g_\xi$ and the connection 
$\xi^{-1}\circ (\pi^*\nabla_B)\circ\xi$. 
This shows that $g_\xi$ is flat. We obtain that a flat 
connection (\ref{isolapos}) gives rise to the pair of a smooth Lorentzian 
metric on $\Sigma\times\R$ and its smooth Levi--Civita connection
\begin{equation}
g_{\xi}=h\circ (\xi\times\xi 
),\:\:\:\:\:\nabla_{B,\xi}=\xi^{-1}\circ (\pi^*\nabla_B)\circ\xi .
\label{konnexio}
\end{equation}
If $\xi$ is not invertible everywhere, the associated metric and 
Levi--Civita connection suffers from singularities. Our construction 
however requires to allow such singular metrics as well hence we will do 
that in what follows. Let us say that two, not necessarily isometric, 
singular metrics are {\it equivalent} if they dreibeins differ only by a 
transformation $\xi\mapsto\xi +(\pi^*\nabla_B)v$. In other words we assign a 
metric to the cohomology class $[\xi ]$ of $\xi$ only. Identifying metrics 
this way has the advantage that although a particular metric (\ref{konnexio}) 
can be singular, within the equivalence class however we can always pass to a 
smooth representative describing an ordinary metric on $\Sigma\times\R$. 

Notice that $g_{\pm\xi}$ are identical metrics (accordingly,
$\nabla_{B,\pm\xi}$ are equivalent).
Define an action of $\Z_2$ on $H^1(\pi^*\nabla_B)$ via $[\xi 
]\mapsto [-\xi ]$. Then the quotient $H^1(\pi^*\nabla_B)/\Z_2$ is 
identified with the space of equivalence classes of flat 
Lorentzian metrics on $\Sigma\times\R$ of the form (\ref{konnexio}).

If $G$ is a Lie group, consider the space
\[{\rm Hom}_0(\pi_1(\Sigma\times\R ), G)/G,\]
where ${\rm Hom}_0$ denotes the discrete 
embeddings of $\pi_1(\Sigma\times\R )\cong\pi_1(\Sigma )$ into $G$. It can 
be identified with one connected component of the space of gauge equivalence 
classes of flat $G$ connections on $\Sigma\times\R$ and has real 
dimension $(2g-2)\dim G$. From our construction 
it is clear that $H^1(\pi^*\nabla_B )$ can be described as the space of 
flat ISO$(2,1)$ connections modulo flat SO$(2,1)$ connections showing its 
real dimension is $h^1=12g-12-(6g-6)=6g-6$. Therefore, putting all 
these things together, we have proved:
\begin{proposition} The first cohomology $H^1(\pi^*\nabla_B)$ of the 
complex (\ref{felhuzas}) admits the following two interpretations. 

First $H^1(\pi^*\nabla_B)$ can be identified with the underlying vector 
space of gauge equivalence classes of those flat ISO$(2,1)$ connections on 
$\widehat{\pi^*E}$, an affine vector bundle with underlying vector bundle 
$\pi^*E$, which are of the form (\ref{isolapos}).

Secondly the quotient $H^1(\pi^*\nabla_B)/\Z_2$ can be identified with the 
space of equivalence classes of inextensible, flat, singular Lorentzian 
metrics on $\Sigma\times\R$ of the form (\ref{konnexio}). 

We have $h^1=6g-6$ for the corresponding Betti number. $\Diamond$
\label{h1}
\end{proposition}
\noindent In light of Proposition \ref{h1} we can 
assign a $6g-6$ dimensional moduli of inequivalent singular Lorentzian 
metrics on $\Sigma\times\R$, 
solutions to the $2+1$ dimensional vacuum Einstein equation, to a 
given flat, irreducible SO$(2,1)$ connection of maximal Euler class. In 
order to achieve a more explicit description of these 
singular metrics, we have to analyze the solutions of the Cartan equation 
(\ref{cartan}) on $\Sigma\times\R$. We can do this by carrying out a 
suitable ISO$(2,1)$ gauge transformation on the connections 
in (\ref{isolapos}).

The splitting $T^*(\Sigma\times\R )\otimes\pi^*E\cong 
(T^*\Sigma\otimes\pi^*E)\oplus (T^*\R\otimes\pi^*E)$ allows us to 
decompose a dreibein
$\xi\in\Omega^1(\Sigma\times\R ,\:\pi^*E)$ as $\xi =\pi^*\xi_t+u_t\dd t$ 
with $\xi_t\in\Omega^1(\Sigma ,\:E)$, 
$u_t\in\Omega^0(\Sigma ,\:E)$ and $t\in\R$. In the obvious temporal 
gauge for $\pi^*\nabla_B$ (see next section), we have $\pi^*\nabla 
'_B=\nabla '_B+\frac{\partial}{\partial t}\dd t$ and then (\ref{cartan}) 
reads as
\[\nabla '_B\xi_t+\left(\frac{\partial\xi_t}{\partial t}+\nabla_B 
u_t\right)\wedge\dd t+\frac{\partial u_t}{\partial t}\:\dd t\wedge\dd 
t=0\]
or simply
\begin{equation}
\nabla '_B\xi_t=0,\:\:\:\:\:\frac{\partial\xi_t}{\partial t}+\nabla_B 
u_t=0
\label{ujcartan}
\end{equation}
over $\Sigma\times\{ t\}$. By fixing a ``Coulomb gauge'' on the 
inhomogeneous part $\xi$ as a next step, we 
can adjust the first equation in (\ref{ujcartan}) into an elliptic one as 
follows. Consider the complex
\begin{equation}
0\rightarrow\Omega^0(\Sigma 
,\:E)\stackrel{\nabla_B}{\longrightarrow
}\Omega^1(\Sigma ,\:E)\stackrel{\nabla '_B}{\longrightarrow
}\Omega^2(\Sigma ,\:E)\rightarrow 0,
\label{komplexus}
\end{equation}
whose pullback is (\ref{felhuzas}). Using the orientation and picking up 
a metric $\gamma$ on $\Sigma$ take the associated elliptic 
complex 
\[0\rightarrow\Omega^1(\Sigma ,\:E)\stackrel{\nabla^*_B\oplus\nabla '_B}
{\longrightarrow}\Omega^0(\Sigma 
,\:E)\oplus\Omega^2(\Sigma, E)\rightarrow 0.\]
We claim that

\begin{proposition} Consider a compact, oriented Riemann surface of 
genus $g>1$ and fix a metric $\gamma$ on it. Let 
$\pi^*\nabla_B$ be an arbitrary flat SO$(2,1)$ connection on 
$\pi^*E\cong T(\Sigma\times\R )$. Then there is a natural vector space 
isomorphism
\[ H^1(\pi^*\nabla_B)\cong{\rm Ker}(\nabla^*_B\oplus\nabla '_B)\]
depending only on the conformal class $[\gamma ]$. That is, for all $[\xi 
]\in H^1(\pi^*\nabla_B)$ there is a unique gauge transformation $\xi 
':=\xi+(\pi^*\nabla_B)v$ with $v\in\Omega^0(\Sigma\times\R , \pi^*E)$ such that 
all solutions of (\ref{cartan}) take the shape $\xi '=\pi^*\eta_{[\xi 
]}$ with 
\[\eta_{[\xi ]}=a_1\eta_1+a_2\eta_2+\dots +a_{6g-6}\eta_{6g-6}\]
where $a_i\in\R$ are constants and $\eta_i\in\Omega^1(\Sigma 
, E)$ with $i=1,\dots,6g-6$ form a fixed basis for the kernel of the 
elliptic operator $\nabla^*_B\oplus\nabla '_B$. That 
is, in this gauge $\xi '$ is independent of time.

Notice that this gauge transformation keeps $\xi$ within its cohomology class 
therefore indeed all solutions of the original Cartan equation (\ref{cartan}) 
over $\Sigma\times\R$ are of this form up to a gauge transformation.
\label{megoldasok}
\end{proposition}
\noindent{\it Proof.} For technical reasons we temporarily put an 
auxiliary Riemannian metric onto $E$ to carry out the calculations in the 
course of this proof. We shall denote by $h^0$, $h^1$ and $h^2$ the 
corresponding Betti numbers of (\ref{komplexus}). The index of this complex 
is equal to  
\[{\rm Index}(\nabla^*_B\oplus\nabla '_B)=-\int\limits_\Sigma 
(3+c_1(E^\C ))\wedge (1+e(T\Sigma ))=-3\cdot (2-2g)=6g-6\]
since $E^\C$ is an PSL$(2,\C )$ bundle consequently its first Chern 
class vanishes (in fact $E^\C$ is a trivial bundle). On the other hand 
${\rm Index}(\nabla^*_B\oplus\nabla 
'_B)=-h^0+h^1-h^2$ and Proposition \ref{h1} shows that 
$h^1=6g-6$ hence we find $h^0=0$ that is, ${\rm
Ker}\nabla_B={\rm Coker}\nabla^*_B=\{0\}$ and $h^2=0$ 
hence ${\rm Coker}\nabla '_B=\{0\}$ showing that actually 
\[{\rm Index}(\nabla^*_B\oplus\nabla '_B)=\dim{\rm Ker}(\nabla^*_B\oplus\nabla 
'_B)=6g-6.\] 
A ``Coulomb'' gauge transformation $\xi '_t:=\xi_t +\nabla_B v_t$ and 
$u'_t:=u_t+\frac{\partial v_t}{\partial t}$
with $v_t=v\vert_{\Sigma\times\{ t\}}\in\Omega^0(\Sigma ,E)$ such 
that $\xi '_t$ satisfies the {\it elliptic} equation
\begin{equation}
(\nabla^*_B\oplus\nabla'_B)\xi '_t=0
\label{dirac}
\end{equation}
exists if and only if $\triangle_Bv_t=-\nabla^*_B\xi_t$ for the gauge 
parameter $v_t$ where $\triangle_B=\nabla^*_B\nabla_B$ is the trace 
Laplacian of $\nabla_B$. 
This equation has solution if $\nabla^*_B\xi_t$ is orthogonal to the 
cokernel of $\triangle_B$ that is, the kernel of $\triangle_B$. However 
${\rm Ker}\triangle_B\cong{\rm Ker}\nabla_B$ which is trivial 
as we have seen hence $\nabla^*_B\xi_t$ is certainly orthogonal to the 
trivial cokernel of $\triangle_B$. Moreover this gauge transformation is 
unique. 

Therefore picking up a fixed basis in the kernel of the elliptic operator 
and observing that $H^1(\pi^*\nabla_B)$ and ${\rm 
Ker}(\nabla^*_B\oplus\nabla '_B)$ are of equal dimensions
we can write all solutions of the first equation of (\ref{ujcartan}) as
\[\xi '_t= f(t)(a_1\eta_1+a_2\eta_2+\dots +a_{6g-6}\eta_{6g-6})\]
with a universal function $f(t)$, independent of $\pi^*\nabla_B$. 
The concrete shape of this function emerges by observing that in this 
gauge, taking into account the second equation of (\ref{ujcartan}) too, we find 
\[\xi '_t\in{\rm Ker}\nabla^*_B,\:\:\:\:\:\frac{\partial\xi 
'_t}{\partial t}\:\bot_{L^2(\Sigma )}\:{\rm Ker}\nabla^*_B\] 
by referring to an $L^2$ scalar product over $\Sigma$. Consequently
\[0=\left\langle \xi '_t\:,\:\frac{\partial\xi '_t}{\partial 
t}\right\rangle_{L^2(\Sigma )} 
=\frac{1}{2}\:\frac{\dd}{\dd t}\Vert\xi '_t\Vert^2_{L^2(\Sigma )}\]
implying $\xi '_t$ is independent of time in this gauge 
therefore we have to set $f(t)=1$. We denote this $\xi '_t$ as 
$\eta_{[\xi ]}$. Furthermore
(\ref{ujcartan}) yields $\nabla_B u'_t=0$ hence $u'_t=0$ by uniqueness as 
claimed. Taking into account the conformal invariance of (\ref{dirac}), 
which is essentially the massless Dirac equation, the result follows. 
$\Diamond$ 
\vspace{0.1in}

\begin{remark}\rm
From the viewpoint of ISO$(2,1)$ gauge theory, we 
interpret this result as the existence of temporal gauge for a connection 
(\ref{isolapos}). Indeed, its $\pi^*\nabla_B$ part is time-independent 
(see next section) as well as the translation $\xi '$ as we have seen.
In light of this proposition a generic representative $\xi\in [\xi ]\in 
H^1(\pi^*\nabla_B)$ looks like
\begin{equation}
\xi =\pi^*(\eta_{[\xi ]}+\nabla_Bv_t)+\frac{\partial v_t}{\partial t}\:\dd t
\label{haromlab}
\end{equation}
and we can suppose that its characteristic part $\eta_{[\xi ]}$, a 
solution of (\ref{dirac}), is always smooth by elliptic regularity. 
Concerning the vector field $v_t$ we only know {\it a priori} that it 
somehow diverges as $t\rightarrow\pm\infty$ because the corresponding 
(singular) metric is inextensible by assumption, but otherwise arbitrary. 
For simplicity we suppose it is smooth (we could relax this regularity).

This shows and we also emphasize that the cohomology class $[\xi ]$ is 
quite immense from a geometric viewpoint: The corresponding
non-isometric flat metrics of the form (\ref{konnexio}) have rich 
asymptotics, depending on the gauge parameter $v$ in (\ref{haromlab}). 

To illustrate this we check some examples. The zero dreibein $0\in 
[0]\in H^1(\pi^*\nabla_B)$ corresponds to the totally degenerate 
``metric'' $g_0=0$ on $\Sigma\times\R$ for an arbitrary flat connection. 

A less trivial example: Let $(\pi^*\nabla_B)v\in [0]\in
H^1(\pi^*\nabla_B)$ be another
representative, a pure gauge still within the zero cohomology class. The
corresponding dreibein arises by taking $\eta_{[0]}=0$ in
(\ref{haromlab}) with a fixed vector field $v_t=v\vert_{\Sigma\times\{ 
t\}}$. Since $E=V\oplus\underline{\R}$ with a vector bundle $V\cong T\Sigma$ 
(non-canonical isomorphism) we have a decomposition 
$\Omega^0(\Sigma ,\:E)\cong\Omega^0(\Sigma
,\:V)\oplus\Omega^0(\Sigma )$. We put simply $v_t:=ta$ with $t>0$ and a 
constant $a=0+1\in\Omega^0(\Sigma ,\:V)\oplus\Omega^0(\Sigma )$. Write
$\gamma$ for the metric whose ``zweibein'' is $\nabla_Ba$ then the 
resulting metric is the incomplete cone metric $-\dd t^2+t^2\gamma$ over
$\Sigma\times\R^+$ as in \cite{fef-gra} (cf. also \cite{and} and 
\cite{wit1}) and in particular $\gamma$ is of constant $-1$ curvature.
\end{remark}

\noindent We provide a further description of  
$H^1(\pi^*\nabla_B)$ which points out its relationship with the 
Teich\-m\"ul\-ler space of $\Sigma$ (also cf. \cite{mon}). Consider a 
cohomology class $[\xi ]\in H^1(\pi^*\nabla_B)$ and 
its unique ``Coulomb'' representative
$\pi^*\eta_{[\xi]}\in [\xi ]$ provided by Proposition \ref{megoldasok}. 
Define
\[V(\nabla_B, [\gamma ]):=\{\eta_{[\xi 
]}\in {\rm Ker}(\nabla^*_B\oplus\nabla '_B)\:\vert\:\mbox{$\eta_{[\xi 
]}$ nowhere vanishes}\} ,\]
the {\it space of regular cohomology classes}. Consider one of its 
connected components $V^+(\nabla_B, [\gamma ])$. Then we assert that

\begin{proposition}
Consider a compact, orientable Riemann surface $\Sigma$ of $g>1$ and 
fix an orientation as well as a metric $\gamma$ and its conformal class 
$[\gamma ]$ on it. Let 
$\nabla_B$ be a smooth, irreducible flat SO$(2,1)$ connection on 
the bundle $E\cong T\Sigma\oplus{\underline\R}$. Consider 
$V^+(\nabla_B,[\gamma ])\subset{\rm Ker}(\nabla^*_B\oplus\nabla '_B)$, 
the regular part of the cohomology group $H^1(\pi^*\nabla_B)$. 
Then there is a natural homeomorphism
\[V^+(\nabla_B, [\gamma ])\cong{\cal T}\]
where ${\cal T}$ denotes the Teichm\"uller space of $\Sigma$ consisting 
of complex structures whose induced orientations agree with the given one 
of $\Sigma$. 
\label{teichmuller}
\end{proposition}
\noindent{\it Proof.} It is very simple. Fix an orientation and a smooth 
metric $\gamma$ on $\Sigma$. Pick up a flat connection $\nabla_B$ on 
the bundle $E\cong T\Sigma\oplus\underline{\R}$. An elliptic 
operator $\nabla^*_B\oplus\nabla '_B$ exists and its kernel depends only 
on $[\gamma ]$ as we have seen. The connection $\nabla_B$ also
induces an ${\rm SO}(2,1)$ structure on $\Sigma$. Take a regular 
solution $\eta_{[\xi]}\in V^+(\nabla_B ,[\gamma ])$ of (\ref{dirac}). 
Since it is parallel i.e., $\nabla'_B\eta_{[\xi ]}=0$ and is nowhere 
vanishing it cuts down the structure group of $T\Sigma$ to 
${\rm SO}(2)\subset {\rm SO}(2,1)$; this together with the 
orientation induces a complex structure $J_{[\xi ]}$ on $\Sigma$. 
The mutually different complex 
structures are enumerated by precisely those representatives which satisfy the 
elliptic equation (\ref{dirac}) i.e., represent nowhere vanishing elements 
in ${\rm Ker}(\nabla^*_B\oplus \nabla '_B)$. $\Diamond$
\vspace{0.1in}

\begin{remark}\rm
We find that $V^+(\nabla_B, [\gamma ])\cong\R^{6g-6}$ 
since the same is true for the Teichm\"uller space (cf. e.g. Corollary 
11.10 in \cite{hit}). 

The cohomology class $[\xi ]\in
H^1(\pi^*\nabla_B)$ consists of very different non-isometric metrics of 
the form (\ref{konnexio}), both singular and regular as we have seen.
Nevertheless we are able to assign a unique boundary conformal class 
at least to each regular cohomology class in a natural way via
Propositions \ref{megoldasok} and \ref{teichmuller}. Indeed, a regular 
cohomology class $[\xi ]$ induces a complex structure $J_{[\xi ]}$ or in 
other words, a conformal class $[\gamma ]_{[\xi ]}$ on $\Sigma$. 
 
However this assignment is less geometric in its nature: The
conformal class to $[\xi ]$ does {\it not} arise by simply restricting a
particular metric within $[\xi ]$ to some boundary at infinity.
\end{remark}

\noindent Summing up we have the following characterization of singular 
Lorentzian metrics; the straightforward verification is left to the reader.
\begin{proposition} 
Let $\pi^*\nabla_B$ be a smooth, irreducible, 
flat SO$(2,1)$ connection 
on the bundle $\pi^*E\cong T(\Sigma\times\R )$ over a compact, oriented 
Riemann surface of genus $g>1$. Consider a cohomology 
class $[\xi ]\in H^1(\pi^*\nabla_B)$ representing an equivalence 
class of singular metrics $g_\xi$ as in (\ref{konnexio}). Also fix a 
conformal class $[\gamma ]$ on $\Sigma$. Then
\begin{itemize}

\item[{\rm (i)}] A generic element $g_\xi$ is an inextensible, smooth, 
symmetric tensor field on $T(\Sigma\times\R )$;

\item[{\rm (ii)}] If $\xi\in [\xi ]$ is moreover invertible 
everywhere then $g_\xi$ is an inextensible, smooth, 
flat, globally hyperbolic Lor\-entz\-ian metric on $\Sigma\times\R$ with 
$\Sigma\times\{ t\}$ representing Cauchy surfaces. The metric $g_\xi$ may 
fail to be complete;

\item[{\rm (iii)}] If in addition $[\xi ]$ represents a regular 
cohomology class in $V^+(\nabla_B,[\gamma ])$ then $[\xi ]$ gives rise 
to a unique conformal class $[\gamma ]_{[\xi ]}$ on $\Sigma$. The regular 
part $V^+(\nabla_B,[\gamma ])$ is parameterized by the Teichm\"uller space 
${\cal T}$ of the oriented Riemann surface. $\Diamond$
\end{itemize}
\label{metrikak}
\end{proposition}

\noindent Before closing this section let us summarize what is known at 
this point. Starting with an irreducible, real solution $[(\nabla_A, 
\Phi )]$ of the SO$(3)$ Hitchin equations over $\Sigma$ which belongs to the 
connected component ${\cal M}_{2g-2}$ of Theorem \ref{hitchintetel}, we have 
found an associated $6g-6$ dimensional moduli of inequivalent singular 
solutions of the vacuum Einstein equation on 
$\Sigma\times\R$ via Proposition \ref{h1}. 
However if we take into account not only the Hitchin pair itself but the 
orientation and the conformal class $[\gamma ]$ on $\Sigma$ as well---which 
are implicitly present---then we can assign to $[(\nabla_A, \Phi )]$ a 
{\it unique} regular cohomology class of singular metrics using 
Propositions \ref{megoldasok} and \ref{teichmuller}. This distinguished 
class of solutions simply arises by picking up {\it that} cohomology 
class $[\xi ]$ on the bulk whose corresponding boundary conformal class 
$[\gamma ]_{[\xi ]}$ in the sense of Propositions \ref{megoldasok} and 
\ref{teichmuller} yields precisely the originally given conformal class 
$[\gamma ]$ on $\Sigma$ regarded as the spacelike past boundary of 
$\Sigma\times\R$.  

In the next section we are looking for the inverse construction.

\section{The inverse construction}

Next we focus our attention to the reverse construction. This will 
turn out to be simple by referring to a powerful theorem of Donaldson. We 
continue to consider compact, orientable Riemann surfaces of genus greater 
than one.

Assume a smooth, flat, probably incomplete singular Lorentzian metric 
$g$ is given on $\Sigma\times\R$ stemming from an irreducible, smooth, flat 
ISO$(2,1)$ connection on an affine bundle 
$\hat{\tilde{E}}$ whose underlying vector bundle is $\tilde{E}\cong 
T(\Sigma\times\R )$. We claim that any flat ISO$(2,1)$ connection is of the 
form (\ref{isolapos}) hence this singular metric and its Levi--Civita 
connection look like (\ref{konnexio}) with a translation $\xi$ and a flat 
connection $\nabla_B$ over the SO$(2,1)$ bundle $E$ on $\Sigma$ whose 
principal bundle is $Q_{2g-2}$ and $\tilde{E}=\pi^*E$. 

Indeed, let $\Gamma$ be a discretely embedded subgroup 
of ${\rm SL}(2,\R )$, isomorphic to $\pi_1(\Sigma )$. Since SL$(2,\R )$ is 
the isometry group of the hyperbolic plane $\HH^2$ we can construct a 
model for $\Sigma$ as the quotient $\HH^2/\Gamma\cong\Sigma$ together with 
the projection $p:\HH^2\rightarrow\Sigma$. We extend this 
to a map $p:\HH^2\times\R\rightarrow\Sigma\times\R$ acting as the 
identity on $\R$.
Then given a flat SO$(2,1)$ connection ${\tilde\nabla}_B$ on 
$\tilde{E}$ its pullback can be written as $p^*\tilde{\nabla}_B =\dd+f^{-1}
\dd f$ with a $\Gamma$-periodic function $f:\HH^2\times\R\rightarrow{\rm SL}
(2,\R )$. We can always gauge away the $\R$-component of the pullback 
connection i.e., we can assume that $f^{-1}\partial_t f=0$ yielding $f$ 
is independent of $t$. Consequently 
in this ``temporal gauge'' the connection $p^*\tilde{\nabla}_B$ hence 
$\tilde{\nabla}_B$ looks like $\pi^*\nabla_B$ 
on $\pi^*E\cong\tilde{E}$ over $\Sigma\times\R$ i.e., gives rise to a 
flat SO$(2,1)$ connection $\nabla_B$ on $E\cong 
T\Sigma\oplus\underline{\R}$. We conclude that a flat ISO$(2,1)$ 
connection is of the form (\ref{isolapos}) with $\nabla_B$ a flat 
connection and $\xi$ a translation. Consequently the associated singular 
metric $g$ possesses the properties summarized in Proposition 
\ref{metrikak} hence we shall denote it as $g_\xi$. 

In particular if $\xi$ is a representative of a regular cohomology class 
in $V^+(\nabla_B,[\gamma ])$ then via Propositions \ref{megoldasok} and 
\ref{teichmuller} it gives rise to a unique boundary conformal 
class $[\gamma ]_{[\xi ]}$ on $\Sigma$, regarded as the spacelike past 
boundary of $\Sigma\times\R$ with induced orientation. Moreover the 
restriction of its Levi--Civita connection to 
$\Sigma$ yields a unique, irreducible flat SO$(2,1)$ connection on $E$. 
Given these data: $[\gamma ]_{[\xi ]}$ and $\nabla_B$ on $\Sigma$ one can 
raise the question whether or not they correspond to a (real) solution of 
the SO$(3)$ Hitchin equations. If yes, then we have $w_2(P)=0$ for the 
corresponding SO$(3)$ principal bundle over $\Sigma$ since the Euler 
class of the underlying SO$(2,1)$ principal bundle of $E=E_{2g-2}$ is 
even.

The question is answered in the affirmative by the following theorem 
\cite{don}:
 
\begin{theorem}{\em (Donaldson, 1987)} Let $P$ be an
SO$(3)$ principal bundle over a compact, oriented Riemann surface
$(\Sigma ,[\gamma ]_{[\xi ]})$. Assume $\nabla_B$ is an irreducible flat 
PSL$(2, \C )$ connection on $P^\C$. Then there exists an PSL$(2, \C )$ gauge
transformation on the complexified bundle $P^\C$ taking the flat 
connection into the form $\nabla_A+\Phi +\Phi^*$ where the pair 
$(\nabla_A, \Phi )$
satisfies the SO$(3)$ Hitchin equations (\ref{hitchin}) with respect to 
the conformal class $[\gamma ]_{\xi ]}$ and the orientation on $\Sigma$. 
$\Diamond$
\label{donaldsontetel}
\end{theorem}
\begin{remark}\rm
If the flat, irreducible PSL$(2,\C )$ connection is real then 
the resulting Hitchin pair is also real in the sense of 
Theorem \ref{hitchintetel} and in particular our real solutions are mapped 
into the ${\cal M}_{2g-2}$ component.
\end{remark}
 
\section{An AdS/CFT-type correspondence}

The time has come to bring all of our findings together. These lead us 
to an AdS/CFT-type correspondence between classical $2+1$ 
dimensional vacuum general relativity on the bulk space-time 
$\Sigma\times\R$ and $2$ dimensional SO$(3)$ Hitchin 
theory---regarded as a classical conformal field 
theory---on the spacelike past boundary $\Sigma$. 

We find the most expressive way to present the duality equivalence by 
formulating it in terms of the corresponding field equations. 
Then we can rephrase it by referring to the solutions themselves. For 
notational simplicity we continue to denote a real Hitchin pair on 
the principal SO$(3)$ bundle $P$ with $w_2(P)=0$ as $(\nabla_A, \Phi )$ 
while $\nabla_A+\Phi +\Phi^*$ is the associated flat connection on the 
SO$(2,1)$ vector bundle $E$ of the SO$(2,1)$ principal bundle 
$Q_{2g-2}$ of Theorem \ref{hitchintetel}. 

\begin{theorem}
Let $(\Sigma ,[\gamma ])$ be an oriented, compact Riemann surface $\Sigma$ 
of genus $g>1$ with a fixed conformal class. Consider 
$[(\nabla_A, \Phi )]\in{\cal M}_{2g-2}$, an irreducible, real
solution of the Hitchin equations on the SO$(3)$ principal bundle $P$ 
over $\Sigma$. Then this pair, consisting of the gauge equivalence class 
of an SO$(3)$ connection $\nabla_A$ and a complex Higgs field $\Phi$, satisfies 
the Hitchin equations over $\Sigma$:
\[\left\{\begin{array}{ll}
                 F_A+[\Phi , \Phi^*] & =0 \\ [2mm]
                 \overline{\partial}_A\Phi & =0. \\ [2mm]
\end{array}\right.\]

Then there is a unique associated pair $[(\pi^*\nabla_B ,\xi )]$ 
consisting of the gauge equivalence class of a flat SO$(2,1)$ connection 
$\nabla_B:=\nabla_A+\Phi +\Phi^*$ on $E\cong 
T\Sigma\oplus\underline{\R}$ and a regular cohomology class 
$[\xi ]\in H^1(\pi^*\nabla_B)$ of a dreibein $\xi 
\in\Omega^1(\Sigma\times\R 
,\:\pi^*E)$ with induced boundary conformal class being precisely 
$[\gamma ]$ such that they satisfy the real vacuum Einstein equation over 
$\Sigma\times\R$ with induced natural orientation:

\[\left\{\begin{array}{ll}
                 \pi^*F_B & =0 \\ [2mm]
                 (\pi^*\nabla '_B)\xi & =0. \\ [2mm]
\end{array}\right.\]

Conversely, given a real solution $[(\pi^*\nabla_B, \xi )]$ of the vacuum 
Einstein equation over the naturally oriented $\Sigma\times\R$ such that
the induced boundary conformal class of $[\xi ]\in H^1(\pi^*\nabla_B)$ 
is precisely $[\gamma ]$, then there 
exists a unique irreducible, real solution $[(\nabla_A,\Phi )]$ of the Hitchin 
equations on the SO$(3)$ principal bundle $P$ over $(\Sigma , 
[\gamma ])$ with induced orientation such that 
$\nabla_A+\Phi +\Phi^*$ is PSL$(2,\C )$ gauge equivalent to $\nabla_B$ 
on $E$. $\Diamond$
\label{ads-cft}
\end{theorem}

\noindent This implies that there is a kind of correspondence between 
certain smooth, real, irreducible solutions of the $2$ dimensional SO$(3)$ 
Hitchin equations and solutions of the $2+1$ dimensional vacuum Einstein 
equation expressed in the more usual form of a metric as follows.

Associated to $[(\nabla_A, \Phi )]\in{\cal M}_{2g-2}$ over the 
oriented $(\Sigma 
,[\gamma ])$ there are singular solutions $g_\xi$ of the 
Lorentzian vacuum Einstein equation on $\Sigma\times\R$ with natural 
induced orientation such that
\begin{itemize}

\item[(i)] The metric and its Levi--Civita connection are of the 
form (\ref{konnexio}) with $\nabla_B=\nabla_A+\Phi +\Phi^*$ and some 
$\xi$. The isometry classes of these singular metrics 
are parameterized by a unique regular cohomology class $[\xi ]\in 
H^1(\pi^*\nabla_B)$;

\item[(ii)] The boundary conformal class, induced by $[\xi ]$ 
in the sense of Propositions \ref{megoldasok} and \ref{teichmuller} on 
$\Sigma$ (regarded as the spacelike past boundary of 
$(\Sigma\times\R ,g_\xi )$), is equal precisely to $[\gamma ]$.
\end{itemize}

\noindent Conversely, given a set of singular solutions $g_\xi$ of the 
Lorentzian vacuum Einstein equation on the naturally oriented $\Sigma\times\R$ 
which together with their Levi--Civita connections are of the 
form (\ref{konnexio}), there is a unique real solution 
$[(\nabla_A, \Phi )]\in{\cal M}_{2g-2}$ of the SO$(3)$ Hitchin equations 
over the spacelike past boundary $(\Sigma ,[\gamma ])$ with induced 
orientation moreover a unique regular cohomology class $[\xi ]\in 
H^1(\pi^*\nabla_B)$ such that  
\begin{itemize}

\item[(i)] The connection $\nabla_A+\Phi +\Phi^*$ is PSL$(2,\C 
)$ gauge equivalent to the Levi--Civita connections of the $g_\xi$'s 
restricted to $\Sigma$ in temporal gauge; 

\item[(ii)] The boundary conformal class induced by this unique regular 
cohomology class $[\xi ]$, in the sense 
of Propositions \ref{megoldasok}
and \ref{teichmuller}, is precisely $[\gamma ]$ on $\Sigma$. 
\end{itemize}
\noindent We can see at this point that this correspondence can be 
interpreted as a sort of {\it generalized} AdS/CFT correspondence between 
these theories.
By ``generalized'' we mean the way of assigning a boundary conformal 
class to a bulk metric: It does {\it not} arise geometrically by taking 
the conformal class of the bulk metric and then restricting one of its 
representatives to the past or future boundary of the bulk. Rather we 
associate the same conformal geometry to 
metrics of probably very different asymptotics, parameterized by a regular 
cohomology class and the conformal geometry arises in an abstract way 
exploiting the conformal properties of a massless Dirac-like
equation as explained in Propositions \ref{megoldasok} and 
\ref{teichmuller}. 

We decided to present the main result in terms of the field 
equations not only because of their impressive form but in this way we 
can also point out that the $2+1$ dimensional vacuum Einstein equation, 
if formulated in terms of a connection and a dreibein, can be viewed as a 
sort of ``decoupled'' version of the SO$(3)$ Hitchin equations: It 
is challenging to view the flat SO$(2,1)$ connection 
$\nabla_B$ as the ``dual'' connection to the non-flat SO$(3)$ connection 
$\nabla_A$ and the dreibein $\xi$ as ``dual'' Higgs field to $\Phi$ and 
vice versa. The straightforward advantage of the Einstein equation over 
the Hitchin equations is that the former is decoupled. Observe that at 
least formally we have no reason to restrict this description to real 
solutions hence this duality can in principle continue to hold for a 
generic complex solution of the Hitchin equations (in the sense that the 
associated flat connection may belong to PSL$(2,\C )$) and for the complex 
dual Higgs field one has $\xi^\C\in\Omega^1_\C (\Sigma\times\R 
,\:\pi^*E^\C )$.

\section{Conclusions}

In this paper we presented a natural classical AdS/CFT-type duality 
between three dimensional Lorentzian vacuum general relativity and two 
dimensional Hitchin conformal field theory. This correspondence might be 
considered as a physical interpretation of at least the real solutions of the 
SO$(3)$ Hitchin equations (cf. the Introduction of \cite{hit}). 

One may try to probe this correspondence beyond the classical 
level by calculating (\ref{particiofv}) in this context. Fix a regular
cohomology class $[\xi ]\in H^1(\pi^*\nabla_B)$ with corresponding 
$[\gamma ]_{[\xi ]}$ on $\Sigma$. Then on the conformal 
side we have the unique data $([\gamma ]_{[\xi ]}, [(\nabla_A,\Phi )])$ 
on $\Sigma$ with a real Hitchin pair while on the gravitational side we 
find Lorentzian metrics 
$g_\xi$ on $\Sigma\times\R$ parametrized by $[\xi ]$. Then the 
partition function of the Hitchin conformal field theory is formally equal 
to
\[Z_{CFT}\left( [\gamma ]_{[\xi ]}, [(\nabla_A, \Phi )]\right) 
=\frac{1}{{\rm Vol}([\xi ])}\int\limits_{[\xi ]}{\rm e}^{{\bf 
i}I(g_\xi )}{\rm D}\xi\]
where the integral is taken over the regular cohomology class $[\xi ]$. A 
generic element is given as in 
(\ref{haromlab}) consequently $[\xi ]\cong\Omega^0(\Sigma\times\R 
,\:\pi^*E)$, an infinite dimensional space. This integral shares some 
similarities with those considered in \cite{moo-nek-sha}. Bearing in 
mind that probably both sides of the above integral expression make no sense 
mathematically, we can calculate it formally as follows. The 
Einstein--Hilbert action on a Lorentzian manifold with 
vanishing cosmological constant and spacelike boundary looks like
\[I(g)=-\frac{1}{16\pi G}\int\limits_{M}s(g) \dd g-\frac{1}{8\pi G}
\int\limits_{\partial M}{\rm tr}k(g)\:\dd (g\vert_{\partial M})\]
with $s(g)$ being the scalar curvature and $k(g)$ the second fundamental 
form of the boundary. In our case $s(g_\xi )=0$ where $\xi$ is 
invertible hence the first term vanishes for regular representatives 
however the second term may not exist for certain asymptotics of 
$v_t$ in (\ref{haromlab}). We can overcome this difficulty if replace 
the action by its holographically renormalizied form $I^{ren}(g_\xi )$ as 
in \cite{wit2} (also cf. \cite{and}); this gives 
simply $I^{ren}(g_\xi )=0$ in our case for all invertible representatives. 
Consequently, by arguing that non-invertible elements of the cohomology 
class form a ``set of measure zero'' we formally find for the 
particular Hitchin pairs in ${\cal M}_{2g-2}$ that
\[Z_{CFT}\left( [\gamma ]_{[\xi ]}, [(\nabla_A, \Phi )]\right)=1.\]

Interesting questions can be raised for future work. For instance, what is 
the physical interpretation of generic complex solutions of the SO$(3)$ 
Hitchin equations? At first sight one can declare without problem that 
they correspond to complex flat metrics on $T^\C (\Sigma\times\R 
)$ but 
this sounds rather unphysical. Taking into account that ${\rm PSL}(2,\C 
)\cong{\rm SO}(3,1)$, the identity component of the {\it four} dimensional 
Lorentz group, one may try to regard the complex solutions as real flat 
metrics on $\Sigma\times\R^2$; however in four dimensions flat metrics do 
not exhaust solutions of the vacuum Einstein equation therefore this 
interpretation would not be ``tight'' enough. Perhaps it is possible 
that a complex solution can be projected somehow to a 
{\it non-flat} real three dimensional connection therefore representing a 
non-vacuum solution or 
rather a solution with non-zero cosmological constant in $2+1$ 
dimensions. The presence of SO$(3,1)$, the de Sitter isometry group, 
suggests this later possibility.

Finally, notice that in fact the whole construction proceeds through a 
complexification phase which cancels out the information 
of the original real group we began with; this was SO$(3)$ in our case 
just because of convenience: Both the 
Hitchin and the Donaldson theorems are formulated with this group. 
However recently new smooth solutions of the SO$(2,1)$ Hitchin equations 
have been discovered \cite{jar-mos} pointing toward the possibility that 
even SO$(2,1)$ Hitchin theory is interesting and can be used to formulate 
a duality if a Donaldson-type theorem could be worked out.
\vspace{0.1in}

\noindent {\bf Acknowledgement.} The author is grateful for the 
stimulating discussions with M. Jardim and R.A. Mosna 
(IMECC-UNICAMP, Brazil) which finally led to the results of this paper
and L. Feh\'er (E\"otv\"os University, Hungary) for some remarks on topology.

The work was supported by CNPq grant No. 150854/2005-6 (Brazil) and OTKA 
grants No. T43242 and No. T046365 (Hungary).

\end{document}